\input amstex
\input amsppt.sty

\def\Aut{\operatorname{Aut}}

\def\Spec{\operatorname{Spec}}
\def\rank{\operatorname{rank}}
\def\phi{\varphi}
\NoBlackBoxes

\topmatter
\title
Proper holomorphic mappings of the spectral unit ball
\endtitle

\author
W\l odzimierz Zwonek
\endauthor

\abstract
We prove an Alexander type theorem for the spectral unit ball $\Omega_n$
showing that there are no non-trivial proper
holomorphic mappings in $\Omega_n$, $n\geq 2$.
\endabstract

\address
Instytut Matematyki, Uniwersytet Jagiello\'nski, Reymonta 4, 30-059 Krak\'ow, Poland
\endaddress
\email
Wlodzimierz.Zwonek\@im.uj.edu.pl
\endemail
\thanks The research was partially supported by the Research Grant No. 1 PO3A 005 28 of the Polish Ministry of Science
and Higher Education.
\endthanks

\thanks
\endthanks

\thanks
keywords: spectral unit ball, proper holomorphic mappings, symmetrized polydisc
\endthanks

\endtopmatter

\document

Let $\Cal M_n$ denote the space of $n\times n$ complex matrices.

In order to avoid
some trivialities and ambiguities we assume in the whole paper that {\bf $n\geq 2$}.

Let $\rho(A):=\max\{|\lambda|:\lambda\in\Spec (A)\}$ be the spectral radius of $A\in\Cal M_n$.
Denote also by $\Spec(A):=\{\lambda\in\Bbb C:\det (A-\lambda \Bbb I_n)=0\}$ the spectrum of $A\in\Cal M_n$,
where the eigenvalues are counted with multiplicities ($\Bbb I_n$ denotes the identity matrix).
We also denote the spectral unit ball by
$\Omega_n:=\{A\in\Cal M_n:\rho(A)<1\}$. Note that $\Omega_n$ is an unbounded pseudoconvex balanced
domain in $\Bbb C^{n^2}$
with the continuous Minkowski functional equal to $\rho$.
For $A\in\Cal M_n$ denote $P_A(\lambda):=\det (\lambda \Bbb I_n-A)=
\lambda^n+\sum\sb{j=1}\sp{n}(-1)^j\sigma_j(A)\lambda^{n-j}$, $A\in\Cal M_n$.
Denote also $\sigma:=(\sigma_1,\ldots,\sigma_n)$.
We put $\Bbb G_n:=\sigma(\Omega_n)$.
The domain $\Bbb G_n$ is called {\it the symmetrized polydisc}. Note that $\sigma\in\Cal O(\Cal M_n,\Bbb G_n)$.
Denote also $\Cal J_n:=\pi_n(\{(\zeta_1,\ldots,\zeta_n):\zeta_j=\zeta_k\text{ for some $j\neq k$}\})$, where
$\pi_{n,j}(\zeta_1,\ldots,\zeta_n):=
\sum\sb{1\leq k_1<\ldots<k_j\leq n}\zeta_{k_1}\cdot\ldots\cdot\zeta_{k_j}$, $\zeta_l\in\Bbb D$, $l=1,\ldots,n$
($\Bbb D$ denotes the unit disc in $\Bbb C$).
Note that $\Bbb G_n\setminus\Cal J_n$ is a domain and $\Bbb G_n\setminus\Cal J_n$ is dense in $\Bbb G_n$.

Note that $\Omega_n=\bigcup\sb{z\in\Bbb G_n}\Cal T_z$, where $\Cal T_z:=\{A\in \Omega_n:\sigma(A)=z\}$, $z\in\Bbb C^n$.
The sets $\Cal T_z$, $z\in\Bbb C^n$ are pairwise disjoint analytic sets.
Note that if the matrix $A\in\Cal T_z$ is
non-degoratory then $A$ is a regular
point of $\Cal T_z$ -- recall that in such a case $\rank\sigma^{\prime}(A)=n$ -- it is the largest possible number.
For definition and basic properties of non-derogatory matrices see \cite{Nik-Tho-Zwo~2007} and references there.
One of possible definitions of a non-derogatory matrix is that different blocks in the Jordan normal form
correspond to different
eigenvalues (or equivalently all eigenspaces  are one-dimensional). We shall deliver some properties of
the sets $\Cal T_z$ (see Lemma 5, Lemma 6 and Corollary 7).
It is also simple to see that $\Cal T_0$ is a cone which contains at least $n^2-n+1$ linearly independent vectors:
for instance the ones consisting of one $1$ lying not on the diagonal (and with other entries equal to $0$)
and the matrix $A=(a_{j,k})_{j,k=1,\ldots,n}$ such that
$a_{1,1}=1$, $a_{1,2}=-1$, $a_{2,1}=1$, $a_{2,2}=-1$ (and with all other entries equal to $0$). Consequently,
we shall see that $0$ is not a regular point of $\Cal T_0$.
On the other hand the sets $\Cal T_{\pi_n(\zeta_1,\ldots,\zeta_n)}$, where the points $\zeta_1,\ldots,\zeta_n\in\Bbb D$
are pairwise different, are submanifolds -- it follows from the fact that in this case all elements
of $\Cal T_{\pi_n(\zeta_1,\ldots,\zeta_n)}$ are non-derogatory.


It is well-known that for a given mapping $F\in\Cal O(\Omega_n,\Omega_n)$ there exists a mapping
$\tilde F\in\Cal O(\Bbb G_n,\Bbb G_n)$ such that $\sigma(F(A))=\tilde F(\sigma(A))$ (see e.g. \cite{Edi-Zwo~2005}).

If $f\in\Cal O(\Bbb D,\Bbb D)$ then one may well-define the following holomorphic
mapping $F_f:\Omega_n\owns A\mapsto f(A):=
\sum\sb{j=0}\sp{\infty}\frac{f^{(j)}(0)}{j!}A^j\in\Omega_n$.
Note that $\tilde F_f(\sigma(A))=\sigma(F_f(A))=\pi_n(f(\lambda_1),\ldots,f(\lambda_n))$,
where $\sigma(A)=\pi_n(\lambda_1,\ldots,\lambda_n)$.
In particular, $F_a\in\Aut(\Omega_n)$ for any $a\in\Aut\Bbb D$. On the other hand the function
$F_{B}$, where $B(\lambda)=\lambda^2$, $\lambda\in\Bbb D$ is a mapping of the form
$\Omega_2\owns A\mapsto A^2\in \Omega_2$, which
is not a proper holomorphic one -- it maps $\Cal T_0$ into $0$.

The structure of the group of automorphisms of $\Omega_n$ has been been studied in several papers
(see e.g. \cite{Ran-Whi~1991}
and \cite{Ros~2003}). However, it is still not understood completely. Let us mention only that $\Aut(\Omega_n)$
is not transitive. Motivated
by the results of the mentioned papers we are going to examine the structure of the class
of proper holomorphic self-mappings
of the spectral unit ball. It turns out that we get an analogue of the theorem of Alexander on proper
holomorphic self mappings of the Euclidean ball in $\Bbb C^n$ stating that there are
no non-trivial proper holomorphic self maps in the unit ball $\Bbb B_n$, $n\geq 2$ (see \cite{Ale~1977}).

In the paper we need some properties of proper holomorphic mappings between complex analytic sets
that could be found in \cite{Chi~1989} and \cite{\L oj~1991}. The book \cite{Rud~1980} may serve as another
reference on proper holomorphic mappings
(mostly between open sets in $\Bbb C^n$).

\proclaim{Theorem 1} Let $F:\Omega_n\mapsto\Omega_n$ be a proper holomorphic mapping, $n\geq 2$.
Then $F$ is an automorphism.
\endproclaim

The following necessary form of proper holomorphic mappings of the spectral ball, which is a simple consequence
of the description of the set of proper holomorphic self-mappings of the symmetrized polydisc, will be crucial
in our considerations and justifies the introducing of the condition \thetag{1} below.

\proclaim{Proposition 2 {\rm (see Theorem 17 in \cite{Edi-Zwo~2005})}} Let $F:\Omega_n\mapsto\Omega_n$
be a proper holomorphic mapping. Then there is a non-constant finite Blaschke product $B$
such that $\sigma(F(A))=\pi_n(B(\zeta_1),\ldots,B(\zeta_n))$, where $A\in\Cal \Omega_n$ and
$\sigma(A)=\pi_n(\zeta_1,\ldots,\zeta_n)$, $\zeta_1,\ldots,\zeta_n\in\Bbb D$.
\endproclaim
In view of Proposition 2 it is natural that we study below the holomorphic mappings
$F:\Omega_n\mapsto\Omega_n$ such that there is a function $f\in\Cal O(\Bbb D,\Bbb D)$ with the property
$$
\sigma(F(A))=\pi_n(f(\zeta_1),\ldots,f(\zeta_n)),\text{ $A\in\Omega_n$ is such that } \sigma(A)=
\pi_n(\zeta_1,\ldots,\zeta_n).\tag{1}
$$


We start with the following lemma.

\proclaim{Lemma 3} Let $F\in\Cal O(\Omega_n,\Omega_n)$ be such that $F(0)=0$ and \thetag{1} is satisfied for
$f\in\Cal O(\Bbb D,\Bbb D)$ (then necessarily $f(0)=0$) with $f^{\prime}(0)\neq 0$.
Then $F^{\prime}(0)$ is a linear isomorphism (of $\Cal M_n$).
\endproclaim
\demo{Proof} Put $\alpha:=f^{\prime}(0)$. Fix $V\in\Cal M_n$. Let $\pi_n(\mu)=\sigma(V)$
for some $\mu=(\mu_1,\ldots,\mu_n)\in\Bbb C^n$.
We first prove that
$$
\frac{1}{\alpha}F^{\prime}(0)(V)\in\sigma(V).\tag{2}
$$
Actually, $\sigma(\zeta V)=\pi_n(\zeta\mu)$,
$\zeta\in\Bbb C$. Consequently,
$$
\sigma(F(\zeta V))=\pi_n(f(\zeta\mu_1),\ldots,f(\zeta\mu_n))
$$
for sufficiently small $\zeta\in\Bbb D$ and then
$$
\sigma\left(\frac{F(\zeta V)}{\zeta}\right)=\pi_n\left(\frac{(f(\zeta\mu_1),\ldots,f(\zeta\mu_n))}{\zeta}
\right).
$$
Passing with $\zeta$ to $0$ we get that  $F^{\prime}(0)(V)\in\Cal T_{\pi_n(\alpha\mu)}$. Therefore,
$\Phi:=\frac{1}{\alpha}F^{\prime}(0):\Cal M_n\mapsto\Cal M_n$ is a linear mapping such that
$$
\Phi(\Cal T_z)\subset\Cal T_z, \;z\in\Bbb C^n.\tag{3}
$$
To finish the proof of the lemma it is sufficient to show that $\Phi$ is a monomorphism.
Suppose that it does not hold. Then there is an $N\in\Cal M_n$, $N\neq 0$
such that $\Phi(N)=0$. Because of \thetag{3} we get that $N\in\Cal T_0$. But then there is an
$M\in\Cal T_0$ such that $N+M\not\in \Cal T_0$.
In particular,
$$
\Cal T_0\not\owns \Phi(N+M)=\Phi(N)+\Phi(M)=\Phi(M)\in\Cal T_0
$$
-- contradiction.



\qed
\enddemo

\proclaim{Lemma 4} Let $F\in\Cal O(\Omega_n,\Omega_n)$
be such that \thetag{1} is satisfied with $f(0)=0$ and $f^{\prime}(0)\neq 0$.
Then $F^{-1}(0)\cap \Cal T_0\subset\{0\}$.
\endproclaim
\demo{Proof} Suppose that there is an $A\in\Cal T_0$, $A\neq 0$ such that $F(A)=0$.
It follows from the Jordan decomposition theorem that there are linearly independent vectors $v_1,v_2\in\Cal M_n$ such that
$A(v_2)=v_1$, $A(v_1)=0$ (at the moment it is essential that $n\geq 2$). Let $(v_1,\ldots,v_n)$ be a
vector base of $\Bbb C^n$.
Define the linear mapping $V:\Bbb C^n\mapsto\Bbb C^n$ (equivalently an element from $\Cal M_n$)
as follows
$V(v_2):=v_1$, $V(v_1):=v_2$ and $V(v_j):=0$, $j=3,\ldots,n$.
Then $(A+\zeta V)^2(v_j)=\zeta(1+\zeta)v_j$, $j=1,2$. Consequently, the properties
of the spectral radius imply that
$$
|\zeta||1+\zeta|\leq\rho((A+\zeta V)^2)\leq\rho^2(A+\zeta V).
$$
For any $\zeta\in\Bbb C$ there are
$\mu_j(\zeta)\in\Bbb C$, $j=1,\ldots,n$ such that $\pi_n(\mu_1(\zeta),\ldots,\mu_n(\zeta))=\sigma(A+\zeta V)$. Then
$\sigma(F(A+\zeta V))=\pi_n(f(\mu_1(\zeta)),\ldots,f(\mu_n(\zeta)))$ for $\zeta\in\Bbb D$ small.
We also know that
$\max\{|\mu_j(\zeta)|:j=1,\ldots,n\}=\rho(A+\zeta V)\geq\sqrt{|\zeta|}\sqrt{|1+\zeta|}$ and
$\rho(A+\zeta V)\to 0$ as $\zeta\to 0$.

Note that
$\rho(\frac{F(A+\zeta V)}{\zeta})\to\rho(F^{\prime}(A)(V))$ as $\zeta\to 0$. But on the other hand
$$
\rho\left(\frac{F(A+\zeta V)}{\zeta}\right)=\max\left\{\frac{|f(\mu_j(\zeta))|}{|\zeta|}:j=1,\ldots,n\right\},
$$
which tends to infinity as $\zeta\to 0$ because $f^{\prime}(0)\neq 0$ -- a contradiction.
\qed
\enddemo



Note that the results proven so far referred to a larger class of mappings than only proper ones. It is possible
that they may have application to the study of more general mappings than only the proper holomorphic ones.

First we show simple results on the geometry of the sets $\Cal T_z$.

\proclaim{Lemma 5} The set of non-derogatory matrices is dense in $\Cal T_z$ for any $z\in\Bbb G_n$.
\endproclaim
\demo{Proof} Fix $z\in\Bbb G_n$. Let $A\in\Cal T_z$. Without loss of generality assume that $A$ is not non-derogatory.
Choose a vector base $\Cal B$ in which $A$ has Jordan normal form.
Let us study two different blocks corresponding to the same $\lambda$ (and the corresponding
vectors from $\Cal B$: $v_1,\ldots,v_k,w_1,\ldots,w_l$, $k,l\geq 1$).
Let $Av_1=\lambda v_1$, $Av_j=\lambda v_j+v_{j-1}$,
$j=2,\ldots,k$, $Aw_1=\lambda w_1$, $Aw_j=\lambda w_j+w_{j-1}$, $j=1,\ldots,l$. For
$\epsilon>0$ define $Bv_1:=\lambda v_1+\epsilon w_l$ and for all other elements of the base $\Cal B$ define $Bv:=Av$,
$v\in\Cal B$, $v\neq v_1$. This easily gives an approximation of $A$ with matrices still in $\Cal T_z$
having one block corresponding
to the eigenvalue $\lambda$ less than in the original matrix. Repeating this procedure for all Jordan blocks
having the same eigenvalues we easily construct a sequence of non-derogatory matrices in $\Cal T_z$ tending to $A$.
\qed
\enddemo

\proclaim{Lemma 6} The set of non-derogatory matrices in $\Cal T_z$ is connected and open in $\Cal T_z$ for any
$z\in\Bbb G_n$.
\endproclaim
\demo{Proof} Fix $z\in\Bbb G_n$. The non-trivial part of the lemma is the connectedness. Let us fix a system of numbers
$(\zeta_1,\ldots,\zeta_n)$ and the sequence of indices $1=k_1<k_2<\ldots<k_{l+1}=n+1$ where
$\zeta_{k_j}=\zeta_{k_j+1}=\ldots=\zeta_{k_{j+1}-1}$, $j=1,\ldots,l$ (and such that no other equalities between
different $\zeta_j$'s hold) and $\pi_n(\zeta_1,\ldots,\zeta_n)=z$. And now for any vector base $(v_1,\ldots,v_n)$ of $\Bbb C^n$
we define the matrix $A$ (more precisely, an element in $\Cal T_z\subset\Cal M_n$) as follows
$A v_{t}=\zeta_jv_t+v_{t-1}$, $j=1,\ldots,l$, $k_j+1\leq t<k_{j+1}$, $A v_{k_j}=\zeta_jv_{k_j}$, $j=1,\ldots,l$.
Note that the above mapping is continuous and its image equals the set of non-derogatory matrices
in $\Cal T_{\pi_n(\zeta_1,\ldots,\zeta_n)}$. This together with the fact that
the set of all vector basis is connected in $\Cal M_n$ completes the proof.
\qed
\enddemo

As a simple corollary of the results on the set of non-derogatory matrices in the sets $\Cal T_z$ we
get the following.

\proclaim{Corollary 7} For any $z\in\Bbb G_n$ the set
$\Cal T_z$ is an analytic irreducible set of codimension $n$.
\endproclaim

At the moment we are ready to move to the proof of our main result.


\demo{Proof of Theorem 1}
First recall that when $F:\Omega_n\mapsto\Omega_n$ is a proper holomorphic mapping then there is a
finite non-constant Blaschke product $B$ such that $\sigma(F(A))=\pi_n(B(\zeta_1),\ldots,B(\zeta_n))$, where
$\sigma(A)=\pi_n(\zeta_1,\ldots,\zeta_n)$. In particular,
$F(\Cal T_{\pi_n(\zeta_1,\ldots,\zeta_n)})\subset\Cal T_{\pi_n(B(\zeta_1),\ldots,B(\zeta_n))}$,
$\zeta_j\in\Bbb D$, $j=1,\ldots,n$. But the properness of $F$ implies even that the equality
$F(\Cal T_{\pi_n(\zeta_1,\ldots,\zeta_n)})=\Cal T_{\pi_n(B(\zeta_1),\ldots,B(\zeta_n))}$, $\zeta_j\in\Bbb D$,
holds -- it is sufficient to note that $\Cal T_z$ is always connected.
Even more,
$F_{|\Cal T_{\pi_n(\zeta_1,\ldots,\zeta_n)}}:\Cal T_{\pi_n(\zeta_1,\ldots,\zeta_n)}\mapsto
\Cal T_{\pi_n(B(\zeta_1),\ldots,B(\zeta_n))}$ is open and proper for any $\zeta_j\in \Bbb D$, $j=1,\ldots,n$.






We claim that for any $\lambda_0\in\Bbb D$ such that $B^{\prime}(\lambda_0)\neq 0$ (note that such points
exist) the function
$$
F_{|\Cal T_{\pi_n(\lambda_0,\ldots,\lambda_0)}} \text{ is injective.}\tag{4}
$$
Actually, making use of the automorphisms of $\Omega_n$ and the properties of Blaschke
products we may assume that $\lambda_0=0$, $B(0)=0$ and $B^{\prime}(0)\neq 0$.
It follows from Lemma 4 that $F^{-1}(0)\cap\Cal T_0=\{0\}$. In particular,
$F(0)=0$.
Now Lemma 3 applies and
we get that $F^{\prime}(0)$ is an isomorphism. Consequently, $F$ is locally invertible near $0$.
Note that there is a neighborhood $\Cal V$ of $0$ such that $\#F^{-1}(A)\cap\Cal T_0=1$ for any $A\in \Cal V\cap\Cal T_0$.
Otherwise there would exist $\Cal T_0\owns A^{\nu},\tilde A^{\nu}$ such that $A^{\nu}\neq \tilde A^{\nu}$ and
$F(A^{\nu})=F(\tilde A^{\nu})\to 0$. But the properness of $F$ implies that (taking if necessary a subsequence)
either both sequences $(A^{\nu}),(\tilde A^{\nu})$ converge to $0$ or at least one of the sequences converges to
a non-zero element $\tilde A$ from $\Cal T_0$ such that $F(\tilde A)=0$. In the first case we contradict
the local invertibility of $F$ near $0$ and in the second case we get two points in $F^{-1}(0)\cap\Cal T_0$
-- a contradiction, too.


Now the analyticity of the set $\{A\in\Cal T_0:\#F^{-1}(A)\cap \Cal T_0=1\}$ (see e.g. \cite{\L oj~1991}, Section V.7.1)
(the mapping
$F_{|\Cal T_0}:\Cal T_0\mapsto\Cal T_0$ is proper and open) and the fact that $\Cal T_0$
is a cone shows that the mapping $F_{|\Cal T_0}:\Cal T_0\mapsto\Cal T_0$ is a one-to-one mapping.

Now we prove the following property.

\item{(5)}
Let $z^{\nu}\to z^0\in\Bbb G_n$, where $z^{\nu}\in\Bbb G_n$ be such that $F_{|\Cal T_{z^{\nu}}}$ is
not injective for any $\nu$ then $F_{|\Cal T_{z^0}}$ is not injective.

Actually, to prove \thetag{5} note that because of the properties of proper holomorphic mappings we
may assume that there are two sequences of non-derogatory matrices $(A^{\nu})$, $(\tilde A^{\nu})$
with $A^{\nu}\neq \tilde A^{\nu}$
lying in $\Cal T_{z^{\nu}}$, $F(A^{\nu})=F(\tilde A^{\nu})$ and
tending to matrices $A,\tilde A\in \Cal T_{z^0}$ such that $A$ is non-derogatory
and $F_{|\Cal T_{z^0}}$ is locally invertible
in $A$. In the case
$A\neq \tilde A$ we are done, so assume that $A=\tilde A$. The local invertibility of $F_{|\Cal T_{z^0}}$
near $A$ implies that there is a neighborhood $\Cal U$ of $A$ such that $F_{|\Cal T_{z^{\nu}}\cap\Cal U}$ is
invertible for $\nu$ large enough, which contradicts the equality $F(A^{\nu})=F(\tilde A^{\nu})$.

We claim that

\item{(6)} for any
$z\in\Bbb G_n$ the mapping $F_{|\Cal T_z}$ is injective.

Put $U:=\{z\in\Bbb G_n:F_{|\Cal T_z}\text{ is injective}\}$.
The fact that $F_{|\Cal T_0}$ is injective shows that $U$ is not empty. The property
\thetag{5} shows that $U$ is open. To see that $U$ is closed in $\Bbb G_n$ take a sequence
$U\owns z^{\nu}\to z\in \Bbb G_n$. Suppose that $z\not \in U$. Then there are different
non-derogatory matrices $A_1,\ldots,A_k\in\Cal T_z$, $C\in\Cal T_w$ with $k\geq 2$ such that $F^{-1}(C)\cap \Cal T_z=
\{A_1,\ldots,A_k\}$. We may choose arbitrarily small
open connected neighborhoods $\Cal U_1,\ldots,\Cal U_k,\Cal V$ of $A_1,\ldots,A_k,C$
such that $\Cal U_l\cap\Cal U_p=\emptyset$ for $l\neq p$,
$\Cal U_j\cap\Cal T_{\tilde z}$ is connected, $\Cal V\cap\Cal T_{\tilde w}$ is connected
for any $\tilde z,\tilde w\in\Bbb G_n$, $j=1,\ldots,k$
and $F^{-1}(\Cal V)\subset\bigcup\sb{j=1}\sp{k}\Cal U_j$.
Consequently, for any $\nu$ there are pairwise disjoint sets
$F(\Cal U_j\cap\Cal T_{z^{\nu}})\cap \Cal V$, $j=1,\ldots,k$
that are open in $\Cal V\cap\Cal T_{w^{\nu}}$ and that are non-empty for $\nu$ large enough. Now the properness of $F$
shows that for $\Cal V$ sufficiently small the sets $F(\Cal U_j\cap \Cal T_{z^{\nu}})$, $j=1,\ldots,k$
cover the whole set $\Cal V\cap\Cal T_{w^{\nu}}$ for $\nu$ large enough; thus contradicting
the connectedness of $\Cal V\cap\Cal T_{w^{\nu}}$.

Since $\Bbb G_n$ is connected we get that $U=\Bbb G_n$, so \thetag{6} is satisfied.




Let $m$ denote the degree of $B$. We claim that
$m=1$. Suppose that $m\geq 2$. Note that taking instead of $F$ the composition of many $F$'s we may assume that
$m\geq n$. There is a point $\zeta_0\in\Bbb D$ such that $B^{-1}(\zeta_0)=\{\zeta_1,\ldots,\zeta_m\}$.
Composing, if necessary, with automorphisms of $\Omega_n$ we may assume that $\zeta_0=0$. Recall that
$\Cal T_{\pi_n(\zeta_1,\ldots,\zeta_n)}$ is an $n^2-n$ dimensional submanifold. Choose
$A\in\Cal T_{\pi_n(\zeta_1,\ldots,\zeta_n)}$ such that $F(A)=0$.
Let $f:=F_{|{\Cal T}_{\pi_n(\zeta_1,\ldots,\zeta_n)}}:\Cal T_{\pi_n(\zeta_1,\ldots,\zeta_n)}\mapsto\Cal T_0$. 
Then $f$ is a holomorphic bijective mapping.
Let us fix a regular point $C$ in $\Cal T_0$. Then the function
$$
\phi:\Bbb C\owns\lambda\mapsto f^{-1}(\lambda C)\in\Cal T_{\pi_n(\zeta_1,\ldots,\zeta_n)}
$$
is holomorphic on $\Bbb C\setminus\{0\}$ (the points $\lambda C$, $\lambda\in\Bbb C\setminus\{0\}$,
are regular in $\Cal T_0$) and continuous at $0$ with $\phi(0)=A$ (use the properness and injectivity of $f$).
Consequently, $\phi$ is holomorphic on $\Bbb C$. Note that $(F\circ\phi)(\lambda)=\lambda C$, $\lambda\in\Bbb C$,
so the tangent space to $\Cal T_{\pi_n(\zeta_1,\ldots,\zeta_n)}$ at $A$ i.e. $T_{A}
(\Cal T_{\pi_n(\zeta_1,\ldots,\zeta_n)})$ is mapped onto
$H:=F^{\prime}(A)(T_A(\Cal T_{\pi_n(\zeta_1,\ldots,\zeta_n)}))$, which contains the vector
$F^{\prime}(A)(\phi^{\prime}(0))=(F\circ\phi)^{\prime}(0)=C$. Consequently, $H$ contains all regular points
of $\Cal T_0$, so it contains the whole $\Cal T_0$, which contains $n^2-n+1$ linearly independent vectors
contradicting the fact that $H$ is at most $n^2-n$ dimensional vector space.

Consequently, we have proven that $\#F^{-1}(C)=1$ for $C\in\Cal \Omega_n$
showing that $F$ is an automorphism.
\qed
\enddemo

\subheading{Acknowledgment} The author wishes to express his gratitude to Witold Jarnicki for fruitful conversations
on the properties of proper holomorphic mappings between analytic sets.

\enddocument